\documentclass[12pt,twoside]{amsart}
\usepackage{amsthm,amsfonts,amssymb,amscd,euscript}

\usepackage[matrix,arrow]{xy}

\author{Florin Ambro} 
\address{DPMMS, CMS\\
University of Cambridge,
Wilberforce Road, Cambridge CB3 0WB, UK.}
\email{f.ambro@dpmms.cam.ac.uk}

\newcommand{\isoto}{{\overset{\sim}{\rightarrow}}}

\newcommand{\Q}{{\mathbb Q}}
\newcommand{\Z}{{\mathbb Z}}
\newcommand{\N}{{\mathbb N}}
\newcommand{\R}{{\mathbb R}}

\newcommand{\calL}{{\mathcal L}}
\newcommand{\calO}{{\mathcal O}}
\newcommand{\calR}{{\mathcal R}}

\newcommand{\frestbar}
 {{\scriptscriptstyle\genfrac{.}{.}{0pt}{10}\shortmid\shortmid}}
\newcommand{\frest}[1]{{}_{\frestbar #1}} 

\newcommand{\bA}{{\mathbf A}}
\newcommand{\bB}{{\mathbf B}}

\newcommand{\bD}{{\mathbf D}}

\newcommand{\bK}{{\mathbf K}}

\newcommand{\bM}{{\mathbf M}}

\newcommand{\codim}{\operatorname{codim}}

\newcommand{\Div}{\operatorname{Div}}

\newcommand{\LCS}{\operatorname{LCS}}

\newcommand{\mult}{\operatorname{mult}}

\newcommand{\rank}{\operatorname{rank}}

\newcommand{\Spec}{\operatorname{Spec}}
\newcommand{\Supp}{\operatorname{Supp}}
\newcommand{\Var}{\operatorname{Var}}

\theoremstyle{plain}
\newtheorem{thm}{Theorem}[section]

\newtheorem{lem}[thm]{Lemma}

\newtheorem{prop}[thm]{Proposition}

\theoremstyle{definition}
\newtheorem{defn}[thm]{Definition}
\newtheorem{defnprop}[thm]{Definition-Proposition}

\newtheorem{exmp}[thm]{Example}
\newtheorem{example}{Example}

\newtheorem{rem}[thm]{Remark}

\newtheorem{ack}{Acknowledgments}   

\theoremstyle{remark}

\newenvironment{sketch}{\begin{proof}[Sketch of proof]}{\end{proof}}

\begin{document}

\bibliographystyle{amsalpha+}
\title[Boundary property]{Shokurov's boundary property}

\begin{abstract} 
For a birational analogue of minimal elliptic 
surfaces $f\colon X\to Y$, the singularities of the 
fibers define a log structure $(Y,B_Y)$ 
in codimension one on $Y$. Via base change, 
we have a log structure $(Y',B_{Y'})$ in codimension
one on $Y'$, for any birational model $Y'$ of $Y$. 
We show that these codimension one log structures glue
to a unique log structure, defined on some birational
model of $Y$ (Shokurov's BP Conjecture). 

We have three applications: inverse of adjunction for
the above mentioned fiber spaces, the invariance of 
Shokurov's FGA-algebras under restriction to exceptional 
lc centers, and a remark on the moduli part of 
parabolic fiber spaces.
\end{abstract}

\maketitle


\setcounter{section}{-1}


\section{Introduction}


\footnotetext[1]{1991 Mathematics Subject Classification. 
Primary: 14J10, 14J17, 14N30. Secondary: 14E30.}

Recall Kodaira's {\em canonical bundle formula} for a 
minimal elliptic surface $f\colon S\to C$ defined over the 
complex number field:
$$
K_S=f^*(K_C+B_C+M_C).
$$
The {\em moduli part} $M_C$ is a $\Q$-divisor such that 
$12M_C$ is integral and 
$\calO_C(12M_C)\simeq J^*\calO_{{\mathbb P}^1}(1)$, 
where $J\colon  C\to {\mathbb P}^1$ is the $J$-invariant 
function. The {\em discriminant} $B_C=\sum_P b_PP$, 
supported by the singular locus of $f$, is computed in 
terms of the local monodromies around the singular fibers 
$S_P$: 
$b_P$ equals 
$\frac{m-1}{m},\frac{1}{2},\frac{1}{6},
\frac{5}{6}, \frac{1}{4},\frac{3}{4},\frac{1}{3},
\frac{2}{3},
$
depending on whether $S_P$ is of type
$
mI_b,I^*_b,II,II^*,III,III^*,IV,IV^*,
$
where $mI_b$ is a multiple fibre of multiplicity $m$,
and $b\ge 0$. Kawamata~\cite{sub2,sba} proposed an 
equivalent definition, which does not require 
classification of the fibers: $1-b_P$ is the 
{\em log canonical threshold} of the log pair 
$(S,S_P)$ in a neighborhood of the fiber $S_P$. 
The minimality may also be removed: the birational 
changes of $S$ over $C$ are controlled by a log pair 
structure on $S$.

The birational analogue consists of data 
$f\colon (X,B)\to Y$, where $f$ is a surjective morphism 
between proper normal varieties and $(X,B)$ is a 
{\em log pair} ($B$ is a $\Q$-Weil divisor such that 
$K+B$ is $\Q$-Cartier), satisfying the following 
properties:

\begin{itemize}
\item[(1)] $(X,B)$ has Kawamata log terminal singularities 
over the generic point of $Y$.
\item[(2)] $\rank f_*\calO_X(\lceil \bA(X,B)\rceil)=1$.
\item[(3)] There exist a positive integer $r$, a rational 
function $\varphi\in k(X)^\times$ and a $\Q$-Cartier divisor 
$D$ on $Y$ such that $$K+B+\frac{1}{r}(\varphi)= f^*D.$$
\end{itemize}

We note here that $B$ need not be effective. However, if
$B$ is effective over the generic point of $Y$, then the
technical assumption (2) is implied by (1). 
The {\em discriminant} on $Y$ of the log 
divisor $K+B$ is the $\Q$-Weil divisor $B_Y:=\sum_P b_P P$,
where $1-b_P$ to be the maximal real number $t$ such 
that the log pair $(X,B+tf^*(P))$ has log canonical 
singularities over the generic point of 
$P$. The sum runs after all codimension one points of $Y$, 
but it has finite support. The {\em moduli part} 
is the unique $\Q$-Weil divisor $M_Y$ on $Y$ satisfying
$$
K+B+\frac{1}{r}(\varphi)= f^*(K_Y+B_Y+M_Y).
$$
If the base space $Y$ is a curve, we recover Kodaira's 
formula (except for effective freeness):

\begin{thm}\label{cv} 
Let $f\colon (X,B)\to Y$ be a fiber space satisfying 
(1),(2),(3) above, with $\dim(Y)=1$. Then the moduli
$\Q$-divisor $M_Y$ is semi-ample, i.e. there exists 
a positive integer $k$ such that $kM_Y$ is Cartier 
and the linear system $|kM_Y|$ is base point free.  
\end{thm}

If the base has dimension at least two, the linear 
system $|kM_Y|$ is expected to have no fixed 
components for $k$ large and divisible. However,
$|kM_Y|$ might have base points, as Prokhorov has
pointed out. To remove the non-determinacy locus
of the expected rational map, we should 
replace $f:(X,B)\to Y$ by a birational base change: if 
$f'\colon X'\to Y'$ is a fiber space induced via 
a birational base change $\sigma\colon Y'\to Y$, we 
have an induced data $f'\colon (X',B_{X'})\to Y'$, 
where $B_{X'}$ is defined by $\mu^*(K+B)=K_{X'}+B_{X'}$:

\[ \xymatrix{
(X,B) \ar[d]_f & (X',B_{X'}) \ar[l]_{\mu} \ar[d]^{f'}\\
Y        & Y' \ar[l]^{\sigma}
} \]

We denote by $B_{Y'}$ and $M_{Y'}$ the discriminant
and moduli part of $f'\colon (X',B_{X'})\to Y'$, 
respectively. The stabilisation of the moduli part
after a certain blow-up is our main result:

\begin{thm}\label{m} 
Let $f\colon (X,B)\to Y$ be a fiber space satisfying 
the above properties (1),(2),(3). Then there exists a proper
birational morphism $Y'\to Y$ with the following properties:
\begin{itemize}
\item[(i)] $K_{Y'}+B_{Y'}$ is a $\Q$-Cartier divisor, 
and $\nu^*(K_{Y'}+B_{Y'})=K_{Y''}+B_{Y''}$ for every
proper birational morphism $\nu\colon Y''\to Y'$. 
\item[(ii)] $M_{Y'}$ is a nef $\Q$-Cartier divisor
and  $\nu^*(M_{Y'})=M_{Y''}$ for every
proper birational morphism $\nu\colon Y''\to Y'$. 
\end{itemize}
\end{thm}

The first part is the positive answer to Shokurov's BP 
Conjecture~\cite[page 92]{Plflips}. Prokhorov and 
Shokurov~\cite{PS} proved (i), by a different method,
in a special case when $X$ is a $3$-fold and $Y$ is a 
surface (they also obtain an explicit description of 
$Y'$). Modulo (i), the second part is a result
of Kawamata~\cite[Theorem 2]{sba}. It is expected that 
$|kM_{Y'}|$ is base point free for a positive integer
$k$. This is known if the generic fibre $F$ of $f$ is a 
curve and $B|_F$ is effective~\cite{koddim,sub2}. 
We have three applications of Theorem~\ref{m}: 

(A) By the definition of the discriminant, $(X,B)$ and 
$(Y,B_Y)$ have ``similar singularities'' ({\em inverse 
of adjunction holds}) only outside a closed subset of 
$Y$ of codimension at least two. We show that inverse 
of adjunction extends to the whole variety $Y$, provided 
that $Y$ is high enough (Theorem~\ref{ia}).

(B) Shokurov~\cite{Plflips} has reduced the existence of 
flips to the finite generatedness of certain (FGA) algebras 
which are asymptotically saturated with respect to 
a Fano variety~\cite[Conjecture 4.39]{Plflips}.
We obtain a descent property for asymptotic saturation 
of algebras (Proposition~\ref{div}). To descend
the numerical assumptions of the FGA/0LP Conjecture,
the semi-ampleness of the moduli $\Q$-divisor $M_{Y'}$
is required (higher codimensional adjunction is expected
to hold for the same reason). However, the nef property 
of the moduli part is sufficient for some applications to 
the Fano case: we show the invariance of (FGA) algebras 
under restriction to exceptional log canonical centers 
(Theorem~\ref{dj}). 

(C) The moduli part of a parabolic fiber space has the
expected Kodaira dimension, provided that the geometric
generic fiber has a good minimal model (Theorem~\ref{k0}).

The proof of Theorem~\ref{m} is based on some of the methods
developed for the proof of Iitaka's Conjecture $C_{n,m}$,
especially ~\cite{koddim} (see ~\cite{Mori} for an 
excellent survey).
The essential ingredient is the universal base 
change for relative canonical divisors, 
in the codimension one semi-stable case. Some of 
the applications in (B) are explained conjecturally 
in ~\cite{Plflips}.

\begin{ack}
I would like to thank Professors Yujiro Kawamata,
Miles Reid and Vyacheslav V. Shokurov for useful 
discussions. This work was supported through a 
European Community Marie Curie Fellowship.
\end{ack}


\section{Preliminary}


A {\em variety} is a reduced and irreducible scheme of 
finite type, defined over an algebraically closed field of 
characteristic zero. An open subset $U$ of a variety $X$
is called {\em big} if $X\setminus U \subset X$ has 
codimension at least two. A {\em contraction} is a proper
morphism $f\colon X\to Y$ such that $\calO_Y=f_*\calO_X$.

Let $\pi\colon X\to S$ be a proper morphism from a normal 
variety $X$, and let $K\in \{\Z,\Q,\R\}$. A $K$-Weil divisor 
is an element of $Z^1(X)\otimes_\Z K$. The {\em round up 
(down)} divisors $\lceil D \rceil$ ($\lfloor D\rfloor$) 
are defined componentwise.
Two $\R$-Weil divisors $D_1, D_2$ are $K$-{\em linearly 
equivalent}, denoted $D_1\sim_K D_2$, if there exist 
$q_i\in K$ and rational functions $\varphi_i\in k(X)^\times$ 
such that $D_1-D_2=\sum_i q_i(\varphi_i)$. An $\R$-Weil divisor 
$D$ is called
\begin{itemize}
\item[(i)] {\em $K$-Cartier} if $D\sim_K 0$ in a neighborhood 
of each point of $X$.
\item[(ii)] {\em relatively nef} if $D$ is $\R$-Cartier and 
$D\cdot C\ge 0$ for every proper curve $C$ contracted by $\pi$.
\item[(iii)] {\em relatively free} if $D$ is a Cartier divisor
and the natural map 
$$
\pi^*\pi_*\calO_X(D)\to \calO_X(D)
$$ 
is surjective.
\item[(iv)] {\em relatively ample} if $\pi$ is a projective
morphism and $D$ belongs to the real cone generated by 
relatively ample Cartier divisors.
\item[(v)] {\em relatively semi-ample} if there exists a 
contraction $\Phi:X\to Y/S$ and a relatively ample 
$\R$-divisor $H$ on $Y$ such that $D\sim_\R \Phi^*H$. If
$D$ is rational, this is equivalent to $mD$ being relatively
free for sufficiently large and divisible positive integers $m$.
\item[(vi)] {\em relatively big} if there exists $C>0$ such that
$\rank \pi_*\calO_X(mD)\ge Cm^d$ for $m$ sufficiently large and 
divisible, where $d$ is the dimension of the generic fibre of
$\pi$.
\end{itemize}

 A divisor $D$ has {\em simple normal crossings} if it 
is reduced and its components are non-singular divisors 
intersecting transversely, in the smooth ambient space $X$.

\begin{defn} (V.V. Shokurov) A {\em $K$-b-divisor} $\bD$ of 
$X$ is a family $\{\bD_{X'}\}_{X'}$ of 
$K$-Weil divisors indexed by all birational models $X'$ 
of $X$, such that $\mu_*(\bD_{X''})=\bD_{X'}$ if
$\mu\colon X''\to X'$ is a birational contraction. 

Equivalently, $\bD=\sum_E \mult_E(\bD) E$ is a $K$-valued 
function on the set of all (geometric) valuations of the field 
of rational functions $k(X)$, having finite support on some 
(hence any) birational model of $X$.
\end{defn}

\begin{example} (1) Let $\omega$ be a top rational differential 
form of $X$. The associated family of divisors 
$\bK=\{(\omega)_{X'}\}_{X'}$ is called the {\em canonical 
b-divisor} of $X$. 

(2) A rational function $\varphi \in k(X)^\times$ defines a 
b-divisor $\overline{(\varphi)}=\{(\varphi)_{X'}\}_{X'}$.

(3) An $\R$-Cartier divisor $D$ on a birational model $X'$ 
of $X$ defines an $\R$-b-divisor $\overline{D}$ such that
$(\overline{D})_{X''}=\mu^*D$ for every birational contraction 
$\mu\colon X''\to X'$. 

(4) For an $\R$-b-divisor $\bD$, the round up (down) 
b-divisors $\lceil \bD \rceil$ ($\lfloor \bD\rfloor$) are 
defined componentwise.
\end{example}

An $\R$-b-divisor $\bD$ is called {\em $K$-b-Cartier} 
if there exists a birational model $X'$ of $X$ such that
$\bD_{X'}$ is $K$-Cartier and $\bD=\overline{\bD_{X'}}$.
In this case, we say that $\bD$ {\em descends to $X'$}.
The {\em relative Kodaira dimension} $\kappa(X/S,\bD)$ of a 
$K$-Cartier b-divisor $\bD$ is defined as the the 
relative Kodaira dimension of $\bD_{X'}$, where 
$X'/S$ is a model where $\bD$ descends.

An $\R$-b-divisor $\bD$ is {\em b-nef}/$S$ 
({\em b-free}/$S$, {\em b-semi-ample}/$S$, {\em b-big}/$S$) 
if there exists a birational contraction $X'\to X$ such that
$\bD=\overline{\bD_{X'}}$, and $\bD_{X'}$ is 
nef (free, semi-ample, big) relative to the induced morphism
$X'\to S$.

To any $\R$-b-divisor $\bD$ of $X$, there is an associated
{\em b-divisorial sheaf} $\calO_X(\bD)$. If $U\subset X$ is 
an open subset, then $\Gamma(U,\calO_X(\bD))$ is the set of 
rational functions $\varphi \in k(X)$ (including $0$) such 
that 
$
\mult_E(\overline{(\varphi)}+\bD)\ge 0
$
for every valuation $E$ with $c_X(E)\cap U\ne \emptyset$.

A {\em log pair} $(X,B)$ is a normal variety $X$ endowed 
with a $\Q$-Weil divisor $B$ such that $K+B$ is $\Q$-Cartier.
A {\em log variety} is a log pair $(X,B)$ such that
$B$ is effective. A {\em relative log pair (variety)}
$(X/S,B)$ consists of a proper morphism $\pi\colon X\to S$ and
a log pair (variety) structure $(X,B)$. The {\em discrepancy
b-divisor} of a log pair $(X,B)$ is
$$
\bA(X,B)=\bK-\overline{K+B}.
$$ 
For a valuation $E$ of $k(X)$, the {\em log discrepancy}
of $E$ with respect to $(X,B)$ is 
$a(E;X,B):=1+\mult_E(\bA(X,B))$. The {\em minimal log 
discrepancy} of $(X,B)$ in a proper closed subset 
$W\subset X$, is 
$$
a(W;X,B):=\inf_{c_X(E)\subseteq W}a(E;X,B).
$$

The {\em lc places} of
$(X,B)$ are valuations $E$ such that $a(E;X,B)=0$, and
their centers on $X$ are called {\em lc centers}. 
We say that $(X,B)$ has {\em log canonical (Kawamata
log terminal) singularities} if $a(E;X,B)\ge 0$ 
($a(E;X,B)> 0$) for every valuation $E$. The
{\em non-klt locus} $\LCS(X,B)$ ({\em non-log canonical locus} 
$(X,B)_{-\infty}$) is the union of all centers $c_X(E)$ of 
valuations $E$ with $a(E;X,B)\le 0$ ($a(E;X,B)< 0$) 
(see~\cite{qlv} for some basic properties). We also denote
$$
\bA^*(X,B)=\bA(X,B)+\sum_{a(E;X,B)=0}E.
$$

A {\em relative generalized log Fano} variety is a relative
log variety $(X/S,B)$ such that $-(K+B)$ is ample/$S$.


\section{The discriminant and moduli b-divisors}


\begin{defn}\label{olp} 
A {\em $K$-trivial fibration} $f\colon (X,B)\to Y$ consists 
of contraction of normal varieties $f\colon X \to Y$ and 
a log pair $(X,B)$, satisfying the following properties:
\begin{itemize}
\item[(1)] $(X,B)$ has Kawamata log terminal 
singularities over the generic point of $Y$.
\item[(2)] $\rank f_*\calO_X(\lceil \bA(X,B)\rceil)=1$.
\item[(3)] There exist a positive integer $r$, a rational 
function $\varphi\in k(X)^\times$ and a $\Q$-Cartier 
divisor $D$ on $Y$ such that 
$$K+B+\frac{1}{r}(\varphi)= f^*D.$$
\end{itemize}
\end{defn}

\begin{rem} The property $\rank f_*\calO_X(\lceil 
\bA(X,B)\rceil)=1$ holds in the following examples:
\begin{itemize}
\item[(a)] $f$ is birational to the Iitaka fibration
of a functional algebra $\calL$ which is 
$\bA(X,B)$-saturated (Lemma~\ref{rkone}).

\item[(b)] $(F,B|_F)$ has Kawamata log terminal 
singularities and $B|_F$ is effective, where $F$ is 
the generic fiber of $f$.

\item[(c)] Let $W$ be the normalization of an exceptional 
log canonical centre of a log variety $(X,B)$, and let 
$h\colon E\to W$ be the unique lc place over $W$.
By adjunction, there exists a $\Q$-divisor $B_E$ 
such that $h\colon (E,B_E)\to W$ is a $K$-trivial
fibration (see ~\cite{qlv}).
\end{itemize}
\end{rem}

 Define $B_Y=\sum_{P\subset Y} b_P P$, where
the sum runs after all prime divisors of $Y$, and  
$$
1-b_P=\sup\{t\in \R; ^\exists U\ni \eta_P, (X,B+tf^*(P)) 
\mbox{ lc sing}/U\}.
$$
The coefficients $b_P$ are well defined, since $(X,B)$
has at most log canonical singularities over the general 
point of $Y$, and each prime divisor is Cartier in a 
neighborhood of its general point. It is easy to see 
that the sum has finite support, so $B_Y$ is a well 
defined $\Q$-Weil divisor on $Y$. By (3), there exists 
a unique $\Q$-Weil divisor $M_Y$ such that the following 
{\em adjunction formula} holds:
$$
K+B+\frac{1}{r}(\varphi)=f^*(K_Y+B_Y+M_Y).
$$

\begin{defn}\cite{thesis}
The $\Q$-Weil divisors $B_Y$ and $M_Y$ are the 
{\em discriminant} and {\em moduli part} of the $K$-trivial 
fibration $f\colon (X,B)\to Y$. Note that $K_Y+B_Y+M_Y$ is 
$\Q$-Cartier.
\end{defn}

The adjunction formula gives a one-to-one correspondence 
between the choices of $M_Y$ and rational functions with 
$\Q$-coefficients $\frac{1}{r}\varphi$ such that
$K_F+B_F+\frac{1}{r}(\varphi|_F)=0$, where $F$ is the 
general fibre of $f$. If $M_Y$ and $M'_Y$ correspond to 
$\frac{1}{r}\varphi$ and $\frac{1}{r}\varphi'$, 
respectively, then there exists a rational function 
$\theta \in k(Y)^\times$ such that 
$\varphi'=\varphi f^*\theta$ and $rM'_Y=(\theta)+rM_Y$.
The smallest possible value 
of $r$ is $b=b(F,B_F)$, uniquely defined by
$$
\{m\in \N; m(K_F+B_F)\sim 0\}=b\N.
$$
Unless otherwise stated, we assume that 
$\frac{1}{r}\varphi$ is fixed, and $r=b(F,B_F)$.

According to the following lemma, $B_Y$ and $M_Y$ are
independent of the choice of a crepant model of $(X,B)$ 
over $Y$:

\begin{lem} Let $\sigma\colon X-\to X'$ be a birational map
defined over $Y$, and let $f'\colon X'\to Y$ be the induced 
morphism. Then there exists a unique $\Q$-Weil divisor 
$B_{X'}$ such that $\sigma\colon (X,B) -\to (X',B_{X'})$ is a 
crepant birational map. Moreover, $(X,B)$ and $({X'},B_{X'})$
induce the same discriminant on $Y$.
\end{lem}

\begin{proof} 
There exists a common normal birational model 
of $X$ and $X'$ which makes following diagram commute:
\[ \xymatrix{
& X'' \ar[dl]_\mu \ar[dr]^{\mu'} &  \\
X \ar[rr]^\sigma \ar[dr]_f & & X' \ar[dl]^{f'} \\
& Y  &
 } \]
Let $K_{X''}+B_{X''}=\mu^*(K+B)$ be the log pullback.
Since $K_{X''}+B_{X''}+\frac{1}{r}(\varphi)= 
{\mu'}^*({f'}^*D)$ and $\mu'$ is birational, we have 
$K_{X''}+B_{X''}={\mu'}^*(K_{X'}+B_{X'})$, where
$B_{X'}:=\mu'_*(B_{X''})$. Therefore there exists a 
crepant log structure on $X'$. The uniqueness of $B_{X'}$ 
is clear.

Finally, note that 
$\mu^*(K+B+tf^*(P))=K_{X''}+B_{X''}+t(f\circ \mu)^*(P)=
{\mu'}^*(K_{X'}+B_{X'}+t{f'}^*(P))$. Therefore the 
thresholds $1-b_P$ induced by $K+B$ and $K_{X'}+B_{X'}$
coincide.
\end{proof}

Let $\sigma\colon Y'\to Y$ be a birational contraction from 
a normal variety $Y'$. Let $X'$ be a resolution of the 
main component of $X\times_Y Y'$ which dominates $Y'$.
The induced morphism $\mu\colon X'\to X$ is birational, and 
let $(X',B_{X'})$ be the crepant log structure on $X'$,
i.e. $\mu^*(K+B)=K_{X'}+B_{X'}$:
\[ \xymatrix{
(X,B) \ar[d]_f & (X',B_{X'}) \ar[l]_{\mu} \ar[d]^{f'}\\
Y        & Y' \ar[l]^{\sigma}
} \]
We say that the $K$-trivial fibration 
$f'\colon (X',B_{X'})\to Y'$ is induced by base change. 
Let $B_{Y'}$ be the discriminant of $K_{X'}+B_{X'}$ on 
$Y'$. Since the definition of the discriminant 
is divisorial and $\sigma$ is an isomorphism over
codimension one points of $Y$, we have 
$B_Y=\sigma_*(B_{Y'})$. This means that there exists 
a unique $\Q$-b-divisor $\bB$ of $Y$ such that 
$\bB_{Y'}$ is the discriminant on $Y'$ of the induced 
fibre space $f'\colon (X',B_{X'})\to Y'$, for 
every birational model $Y'$ of $Y$. We call $\bB$ the 
{\em discriminant $\Q$-b-divisor} induced by $(X,B)$ on 
the birational class of $Y$.
Similarly, if we fix the $\Q$-rational function 
$\frac{1}{r}\varphi$, there exists a unique 
{\em $\Q$-b-divisor} $\bM$ of $Y$ such that
$$
K_{X'}+B_{X'}+\frac{1}{r}(\varphi)=
f^*(K_{Y'}+\bB_{Y'}+\bM_{Y'})
$$
for every $K$-trivial fibration $f'\colon (X',B_{X'})\to Y'$
induced by base change on a birational model $Y'$
of $Y$. We call $\bM$ the {\em moduli $\Q$-b-divisor} 
of $Y$, induced by $f\colon (X,B)\to Y$.
We restate Theorem~\ref{m} in terms of b-divisors:

\begin{thm}\label{main} Let $f\colon (X,B)\to Y$ be a
$K$-trivial fibration, and let $\pi\colon Y\to S$ be a 
proper morphism. Let $\bB$ and $\bM$ be the 
induced discriminant and moduli $\Q$-b-divisors of $Y$. 
Then
\begin{itemize}
\item[(1)] $\bK+\bB$ is $\Q$-b-Cartier.
\item[(2)] $\bM$ is b-nef/$S$.
\end{itemize}
\end{thm}
 
We expect Theorem~\ref{main} to hold if we 
allow $\R$-boundaries and $\R$-linear equivalence 
instead of $\Q$-boundaries and $\Q$-linear equivalence
in Definition~\ref{olp}, or if $(X,B)$ has log canonical
singularities over the generic point of $Y$. 
In the latter case, the assumption 
$\rank f_*\calO_X(\lceil 
\bA(X,B)\rceil)=1$ should be replaced by
$
\rank f_*\calO_X(\lceil \bA^*(X,B)\rceil)=1.
$


\section{Inverse of adjunction}


Let $f\colon (X,B)\to Y$ be a $K$-trivial fibration.
The $\Q$-b-divisor of $Y$
$$
\bA_{div}:=-\bB
$$ 
is called the {\em divisorial discrepancy b-divisor}
~\cite[page 92]{Plflips}.
Theorem~\ref{main}(1) is equivalent to the following 
property: there exists a birational model $Y'$ of $Y$ such
that 
$
\bA_{div}=\bA(Y'',\bB_{Y''})
$
for every birational model $Y''$ which dominates $Y'$.
As a corollary, inverse of adjunction holds for 
$f\colon (X,B)\to (Y,\bB_Y)$, after a sufficiently high 
birational base change:

\begin{thm} (Inverse of adjunction)\label{ia}
Let $f\colon (X,B)\to Y$ be a $K$-trivial fibration
such that $\bA_{div}=\bA(Y,\bB_Y)$.
Then there exists a positive integer $N$ such that 
$$
\frac{1}{N}a(f^{-1}(Z);X,B)\le a(Z;Y,B_Y)\le 
a(f^{-1}(Z);X,B).
$$
for every closed subset $Z\subset Y$, where
$a(Z;Y,B_Y)$ and $a(f^{-1}(Z);X,B)$ are the minimal 
log discrepancies of $(Y,B_Y)$ in $Z$, and
$(X,B)$ in $f^{-1}(Z)$ respectively.

In particular, $(Y,B_Y)$ has Kawamata log terminal 
(log canonical) singularities in a neighborhood of 
a point $y\in Y$ if and only if $(X,B)$ has Kawamata 
log terminal (log canonical) singularities in a 
neighborhood of $f^{-1}(y)$.
\end{thm}

\begin{proof} The assumption $\bA_{div}=\bA(Y,\bB_Y)$ 
means that the 
Base Change Conjecture~\cite[Section 3]{thesis} holds for 
$f\colon (X,B)\to Y$. The claim is proved in 
~\cite[Proposition 3.4]{thesis}, but with $N$ depending 
on $Z$. The possible values for minimal log 
discrepancies of a fixed log pair are 
finite~\cite[Theorem 2.3]{mld}, hence a maximal value 
$N=\max_{Z\subset Y} N(Z)$ exists.
\end{proof}

\begin{lem}\label{a1} 
Let $f\colon (X,B)\to Y$ be a $K$-trivial fibration
such that $\bA_{div}=\bA(Y,\bB_Y)$. Assume moreover
that $X,Y$ are non-singular varieties, and the divisors
$B,\bB_Y$ have simple normal crossings support.
Then $f_*\calO_X(\lceil -B\rceil)=\calO_Y(\lceil 
-\bB_Y\rceil)$.
\end{lem}

\begin{proof} By (i), $\calO_X(\lceil \bA(X,B)\rceil)=
\calO_X(\lceil -B\rceil)$. Since $B$ has Kawamata log terminal
singularities over the generic point of $Y$, we have a natural
inclusion 
$$
\calO_Y|_V \subseteq f_*\calO_X(\lceil -B\rceil)|_V
$$
for some open subset $V\subset Y$. Since
$\rank f_*\calO_X(\lceil \bA(X,B)\rceil)=1$,
the above inclusion is an equality, after possibly 
shrinking $V$. Thus we identify
$f_*\calO_X(\lceil -B\rceil)$ with a subsheaf of the
constant sheaf $k(Y)$.
We first show that 
$f_*\calO_X(\lceil -B\rceil)\subseteq \calO_Y
(\lceil -\bB_Y \rceil)$. Let $\varphi$ be a 
rational function of $Y$ such that 
$(f^*\varphi)+\lceil -B \rceil\ge 0$, and let $P$
be a prime divisor of $Y$. We may replace $X$ by some 
resolution, so that there exists a prime divisor $Q$ of 
$X$ such that $f(Q)=P$ and 
$$
1-\mult_P(\bB_Y)=\frac{1-\mult_Q(B)}{m_{Q/P}},
$$
where $m_{Q/P}$ is the multiplicity of $f^*(P)$ at 
$Q$. By assumption, we have $\mult_Q(f^*\varphi)+1-
\mult_Q(B)>0$. But $\mult_Q(f^*\varphi)=m_{Q/P}
\cdot \mult_P(\varphi)$, hence
$\mult_P(\varphi)+1-\mult_P(\bB_Y)>0$. Therefore
$(\varphi)+\lceil -\bB_Y\rceil$ is effective at $P$.

Conversely, assume $(\varphi)+\lceil -\bB_Y\rceil$ is
effective, and fix a prime divisor $Q$ of $X$. There
exists a birational base change
\[ \xymatrix{
(X,B) \ar[d]_f & (X',B_{X'}) \ar[l]_{\mu} \ar[d]^{f'}\\
Y        & Y' \ar[l]^{\sigma}
} \]
such that $P:=f(Q)$ is a prime divisor of $Y'$. 
We have $\sigma^*(K_Y+\bB_Y)=K_{Y'}+\bB_{Y'}$ by 
$\bA_{div}=\bA(Y,\bB_Y)$. Furthermore, the simple
normal crossings assumption implies
$
\sigma_*\calO_{Y'}(\lceil -\bB_{Y'}\rceil)=
\calO_Y(\lceil -\bB_Y\rceil).
$
Therefore $(\varphi)+\lceil -\bB_{Y'}\rceil \ge 0$, hence
$\mult_P(\varphi)+1-\mult_P(\bB_{Y'})>0$. Since
$$
1-\mult_P(\bB_Y) \le \frac{1-\mult_Q(B_{X'})}{m_{Q/P}},
$$
we infer $\mult_Q(f^*\varphi)+1-\mult_Q(B_{X'})>0$, i.e.
$(f^*\varphi)+\lceil -B\rceil $ is effective at $Q$.
\end{proof}

\begin{rem}\label{a2} 
Let $f\colon (X,B)\to Y$ be a $K$-trivial fibration.
If $L$ is a $\Q$-Cartier divisor on $Y$, let $B':=B+f^*L$. 
Then $f\colon (X,B')\to Y$ is a $K$-trivial fibration, 
with moduli $\Q$-b-divisor $\bM'=\bM$, and discriminant 
$\Q$-b-divisor $\bB'=\bB+\overline{L}$.
\end{rem}


\section{Covering tricks and base change}


\begin{thm} \label{keycover}
\cite{abelvar} Let $X$ be 
a non-singular quasi-projective variety endowed with a 
divisor $D$ with simple normal crossings singularities,
and let $N$ be a positive integer. Then there exists a
finite Galois covering $\tau\colon \tilde{X}\to X$ 
satisfying the following conditions:

\begin{itemize}
\item[(1)] $\tilde{X}$ is a non-singular quasi-projective
variety, and there exists a simple normal crossings
divisor $\Sigma_X$ such that $\tau$ is \'etale over
$X\setminus \Sigma_X$, and $\tau^{-1}(\Sigma_X)$ is a 
divisor with simple normal crossings. 

\item[(2)] The ramification indices of $\tau$ over the
prime components of $D$ are divisible by $N$.

\end{itemize}
\end{thm}

\begin{sketch} We may assume that $X$ is projective
(by Hironaka's resolution of singularities, we can compactify 
to complement of snc in projective, construct the cover, and 
then restrict back to the original variety).
Let $A$ be a very ample divisor such that $NA-D_i$ is 
very ample for each component $D_i$ of $D$. Let $n=\dim(X)$. 
There exists 
$H^{(i)}_1,\ldots,H^{(i)}_n \in |NA-D_i|$ for every $D_i$,
such that $\Sigma_X:=D+\sum_{i,j} H^{(i)}_j$ is a divisor 
with simple normal crossings.
Let $X=\cup U_\alpha$ be an affine cover, and let
$D_i+H^{(i)}_j=(\varphi^{(i)}_{j\alpha})$ on $U_\alpha$.
The field extension 
$L:=k(X)[(\varphi^{(i)}_{j\alpha})^{\frac{1}{N}}; i,j]$
is independent of the choice of $\alpha$. 
Let $\tilde{X}$ be the normalization on $X$ in $L$.
Then $\tilde{X}$ is non-singular and
$\tau$ is a Kummer cover, \'etale outside $\Sigma_X$, 
and $\tau^{-1}(\Sigma_X)$ has simple normal crossings.
\end{sketch}

\begin{rem}\label{comp} In the above notations, assume that 
$\varrho\colon Y\to X$ is a surjective morphism from a 
non-singular quasi-projective variety $Y$ such that
$\varrho^{-1}(D)$ has simple normal crossings. Then 
we may assume that $\tau\colon \tilde{X}\to X$ fits into 
a commutative diagram
\[ \xymatrix{
\tilde{X} \ar[d]_\tau & \tilde{Y} \ar[l]_g \ar[d]_{\nu} \\
X        & Y \ar[l]_{\varrho}
} \]
satisfying the following properties:
\begin{itemize}
\item[(1)] $\nu$ is a finite covering and $g$ is a 
projective morphism.
\item[(2)] $\tilde{Y}$ is non-singular quasi-projective.
\item[(3)] There exists a simple normal crossings
divisor $\Sigma_Y$ such that $\nu$ is \'etale over 
$Y\setminus \Sigma_Y$,
$\nu^{-1}(\Sigma_Y)$ has simple normal crossings,
and $\varrho^{-1}(\Sigma_X)\subseteq \Sigma_Y$.
\end{itemize}
\end{rem}

\begin{proof} In the proof Theorem~\ref{keycover}, we 
may choose the divisors $H^{(i)}_j$ so that
$\varrho^{-1}(D+\sum_{i,j} H^{(i)}_j)$ is a divisor with 
simple normal crossings on $Y$. Let
$\tau_1\colon \bar{Y}\to Y$ be the normalization of the main
component of the pull back of $\tau$ to $Y$. 
\[ \xymatrix{
\tilde{X} \ar[d]_\tau & \bar{Y} \ar[l] \ar[d]_{\tau_1} & 
\tilde{Y} \ar[l] \ar[d] \\
X        & Y \ar[l]_{\varrho}  & Y_1 \ar[l]_\pi
} \]
Then $\tau_1$ is a finite cover whose ramification locus is
contained in the simple normal crossings divisor 
$\varrho^{-1}(\Sigma_X)$. Let $N'$ be
the least common multiple of its ramification indices,
and construct by Theorem~\ref{keycover} a finite cover 
$\pi\colon Y_1\to Y$ with respect to $\varrho^{-1}(\Sigma_X)$ 
and $N'$. 
Let $\tilde{Y}/Y_1$ be the normalization of the main component
of the pull back of $\tau_1$ to $Y_1$. 
The induced map $\nu\colon \tilde{Y}\to Y$ is a finite cover.
By construction, $\tilde{Y}/Y_1$ is \'etale, hence $\tilde{Y}$ 
is non-singular. There exists a simple normal crossings divisor
$\Sigma_Y$ containing $\varrho^{-1}(\Sigma_X)$ such that 
$\pi$ is \'etale over $Y\setminus \Sigma_Y$ and 
$\pi^{-1}(\Sigma_Y)$ has simple normal crossings. 
Therefore $\nu^{-1}(\Sigma_Y)$ has simple normal crossings.
\end{proof}

\begin{thm}\label{ss} \cite{KKMS,wp}
(Semi-stable reduction in codimension one) 
Let $f\colon X\to Y$ be a surjective morphism of non-singular 
varieties.
Assume $\Sigma_X,\Sigma_Y$ are simple normal crossings
divisors on $X$ and $Y$ respectively, such that 
$f^{-1}(\Sigma_Y)\subseteq \Sigma_X$ and $f$ is smooth 
over $Y\setminus \Sigma_Y$. 
Then there exists a positive integer $N$ such that the 
following hold: 

Let $\pi\colon Y'\to Y$ be a finite covering from a nonsingular 
variety $Y'$ such that $\Sigma_{Y'}:=\pi^{-1}(\Sigma_Y)$ 
has simple normal crossings and $N$ divides the 
ramification indices of $\pi$ over the prime components
of $\Sigma_{Y'}$. Then there exists a commutative diagram
\[ \xymatrix{
X \ar[d]_f & X\times_Y Y' \ar[l] \ar[d] & X' \ar[l]_p 
\ar[dl]^{f'}\\
Y        & Y' \ar[l]^{\pi} & 
} \]
with the following properties:
\begin{itemize}
\item[(a)] $X'$ is non-singular and 
$\Sigma_{X'}:={\pi'}^{-1}(\Sigma_X)$ has
simple normal crossings, where $\pi'\colon X'\to X$ is the 
induced projective morphism.

\item[(b)] $p$ is projective, and is an isomorphism above
$Y'\setminus \Sigma_{Y'}$. In particular, $f'$ is smooth over
$Y'\setminus \Sigma_{Y'}$. 

\item[(c)] $f'$ is {\em semi-stable in codimension one}:
the fibers over (generic) codimension one points
of $Y'$ have simple normal crossings singularities.
\end{itemize}
\end{thm}

\begin{sketch} 
Let $f^*(\Sigma_X)=\sum n_i E_i$, and let
$N$ be the least common multiple of the $n_i$'s corresponding
to components $E_i$ which dominate some component of $\Sigma_Y$. 
Consider a finite base change $Y'\to Y$ as above. 
Then ~\cite[IV]{KKMS}
shows that over the generic point of each prime component $Q$
of $\Sigma_{Y'}$, $X\times_Y Y'$ admits a resolution with the 
desired properties. Therefore there exists a 
closed subscheme $B \subset X\times_Y Y'$, supported
over $\Sigma_{Y'}$, and a closed subset $Z\subset \Sigma_{Y'}$ 
with $\codim(Z,Y')\ge 2$, such that the blow-up of 
$X\times_Y Y'$ in $B$ has the desired properties over
$Y'\setminus Z$. Then we may take $X'$ to be any resolution 
of the blow-up, which is an isomorphism outside its singular 
locus, and such that (a) holds.
\end{sketch}

\begin{thm}\cite{Fjo,abelvar}\label{ubc}
Let $f\colon X\to Y$ be a projective morphism of non-singular 
algebraic varieties. Assume $f$ is semi-stable in codimension 
one, and there exists a simple normal crossings divisor 
$\Sigma_Y$ such that $f$ is smooth over $Y\setminus \Sigma_Y$.
Then the following properties hold:
\begin{itemize}
\item[(1)] $f_*\omega_{X/Y}$ is a locally free sheaf on $Y$.

\item[(2)] $f_*\omega_{X/Y}$ is {\em semi-positive}:
let $\nu\colon  C \to Y$ be a proper morphism from
a non-singular projective curve $C$, and let $\calL$ be an
invertible quotient of $\nu^*(f_*\omega_{X/Y})$. Then
$\deg(\calL)\ge 0$.

\item[(3)] Let $\varrho\colon  Y'\to Y$ be a projective morphism
from a non-singular variety $Y'$ such that 
$\varrho^{-1}(\Sigma_Y)$ is a simple normal crossings divisor.
Let $X'\to (X\times_Y Y')_{main}$ be a resolution of the
component of $X\times_Y Y'$ which
dominates $Y'$, and let $h\colon X'\to Y'$ be the induced fibre 
space:
\[ 
\xymatrix{
 X \ar[d]_f  & X' \ar[l]    \ar[d]_{f'}   \\
 Y  & Y' \ar[l]_{\varrho}
} \]
Then there exists a natural isomorphism
$\varrho^*(f_*\omega_{X/Y})\isoto f'_*\omega_{X'/Y'}$
which extends the base change isomorphism over 
$Y\setminus \Sigma_Y$.
\end{itemize}
\end{thm}

\begin{sketch} By the Lefschetz principle and flat base 
change, we may assume $k=\mathbb C$. Let 
$Y_0=Y\setminus \Sigma_Y, X_0=f^{-1}(Y_0)$, and let 
$d=\dim(X/Y)$. The locally free sheaf 
$H_0:=R^df_*{\mathbb \Q}_{X_0}
\otimes_{\Q_{Y_0}} \calO_{Y_0}$ is endowed with the 
integrable Gauss-Manin connection and is the underlying
space of a variation of Hodge structure of weight $d$
on $Y_0$, with $F^dH_0=f_*\omega_{X_0/Y_0}$. 
Since $f$ is semi-stable in codimension one, $H_0$ has 
unipotent local monodromies around the components of 
$\Sigma_Y$. 
Let $H$ be the {\em canonical extension}~\cite{Eqdif} of 
$H_0$. By Schmid's asymptotic behaviour
of variations of Hodge structure, the natural inclusion
$$
f_*\omega_{X/Y} \to j_*(F^dH_0)\cap H
$$
is an isomorphism and $f_*\omega_{X/Y}$ is locally free
~\cite{abelvar}. The semi-positivity follows from 
unipotence and Griffiths' semi-positivity of the curvature 
of the last piece of a variation of Hodge structure
~\cite{Fjo,abelvar}. 

For base change, the sheaf $f'_*\omega_{X'/Y'}$ is 
independent of birational changes in $X'$ over $Y'$.
Thus we may assume that 
$X' \to X\times_Y Y'$ is an isomorphism above 
$Y'\setminus \Sigma_{Y'}$, where 
$\Sigma_{Y'}=\varrho^{-1}(\Sigma_Y)$.
Let $H'_0$ be the variation of Hodge structure on
$Y'\setminus \Sigma_{Y'}$ induced by $f'$, and let 
$H'$ be its canonical extension to $Y'$. Since $H_0$ has 
unipotent local monodromies around the components of
$\Sigma_Y$, the canonical extension is compatible with 
base change~\cite[Proposition 1]{koddim}:
$$
H' \isoto \varrho^*H.
$$
This isomorphism preserves the extensions of the 
Hodge filtration, hence it induces an isomorphism
$
\varrho^*(f_*\omega_{X/Y})\isoto f'_*\omega_{X'/Y'}.
$
\end{sketch}

\begin{thm}\cite{Fj, cvbase}\label{fk} 
Let $f\colon X\to Y$ be a contraction from a non-singular 
projective variety $X$ to a projective curve $Y$, and let 
$E$ be a quotient locally free sheaf of $f_*\omega_{X/Y}$. 
If $\deg(\det(E))=0$, then $\det(E)^{\otimes m}\simeq 
\calO_Y$ for some positive integer $m$.
\end{thm}

\begin{sketch} By ~\cite{Fj}, $E$ is a local system which 
is a direct summand of $f_*\omega_{X/Y}$. Since
$E|_{Y_0}$ is a local subsystem of the variation of Hodge 
structure $H_0$, $\det(E|_{Y_0})^{\otimes m}\simeq 
\calO_{Y_0}$ for some positive integer $m$~\cite{Eqdif}.
By flatness, $\det(E)^{\otimes m}\simeq \calO_Y$.
\end{sketch}


\section{An auxiliary relative $0$-log pair}


We prove Theorems~\ref{main}, ~\ref{cv} in this section. 
The following finite base change formula is essential:

\begin{lem}\cite[Theorem 3.2]{thesis} \label{kfbc} 
Consider a commutative diagram of normal varieties
\[ \xymatrix{
(X,B) \ar[d]_f & (X',B_{X'}) \ar[l]_{\nu} \ar[d]^{f'}\\
Y        & Y' \ar[l]^{\tau}
} \]
with the following properties:
\begin{itemize}
\item[(1)] $(X,B)$ is a log pair with log canonical 
singularities over the generic point of $Y$.
\item[(2)] $\tau$ is a finite morphism and $\nu$ is 
generically finite, $f,f'$ are proper surjective.
\item[(3)] $\nu^*(K+B)=K_{X'}+B_{X'}$.
\end{itemize}

Let $B_Y$ and $B_{Y'}$ be the discriminants of 
$K+B$ and $K_{X'}+B_{X'}$ on $Y$ and $Y'$ respectively.
Then $\tau^*(K_Y+B_Y)=K_{Y'}+B_{Y'}$ (pull back of 
$\Q$-Weil divisors under a finite morphism).
\end{lem}

The category of $K$-trivial fibrations is closed
under generically finite base changes. In order to 
{\em normalize} the discriminant $B_Y$ and the 
moduli part $M_Y$, we have to replace the generic 
fibre of $X/Y$ by a generically finite cover. 
The property 
$\rank f_*\calO_X(\lceil \bA(X,B)\rceil)=1$ is not
invariant under this operation, thus we will 
consider an auxiliary fibre space (cf.~\cite{koddim, 
Mori}).
Throughout this section we consider the following 
{\em set-up}:
\[ \xymatrix{
(X,B) \ar[d]_f & \tilde{X} \ar[l]_\pi 
\ar[dl]_{\tilde{f}} & (V,B_V) \ar[l]_d \ar[dll]^h  \\
Y &  & & 
} \]
$f\colon (X,B)\to Y$ is a $K$-trivial fibration, 
$b=b(F,B_F)$ and
$$
K+B+\frac{1}{b}(\varphi)=f^*(K_Y+B_Y+M_Y)
$$
$\pi\colon \tilde{X}\to X$ is the normalization
of $X$ in $k(X)(\varphi^\frac{1}{b})$ and 
$d\colon V -\to \tilde{X}$ is a birational map from a 
non-singular variety $V$.
The induced rational map $g\colon  V -\to X$ is generically 
finite, so there exists a unique log structure $(V,B_V)$ 
such that $g\colon (V,B_V)-\to (X,B)$ is crepant.
We assume the following properties hold:

\begin{itemize}
\item[(i)] $X,V,Y$ are non-singular quasi-projective varieties
endowed with simple normal crossings divisors $\Sigma_X,\Sigma_V,
\Sigma_Y$ on $X$, $V$ and $Y$ respectively.

\item[(ii)] $f$ and $h$ are projective morphisms.

\item[(iii)] $f$ and $h$ are smooth over $Y\setminus \Sigma_Y$, 
and $\Sigma_X^h/Y$ and $\Sigma_V^h/Y$ have relative simple 
normal crossings over $Y\setminus \Sigma_Y$.

\item[(iv)] $f^{-1}(\Sigma_Y)\subseteq \Sigma_X$,
$f(\Sigma_X^v)\subseteq \Sigma_Y$ and
$h^{-1}(\Sigma_Y)\subseteq \Sigma_V$,
$h(\Sigma_V^v)\subseteq \Sigma_Y$.

\item[(v)] $B$, $B_V$ and $B_Y,M_Y$ are supported by 
$\Sigma_X$, $\Sigma_V$ and $\Sigma_Y$, respectively.

\end{itemize}

In this context, the properties (1) and (2) in the 
definition of the $K$-trivial fibration $f:(X,B)\to Y$
are equivalent to the following properties:  
$\lceil -B_F\rceil$ is an effective divisor and
$\dim_k H^0(F,\lceil -B_F \rceil)=1$.

\begin{lem}\label{ex} The following properties hold for
the above set-up:
\begin{itemize}
\item[(1)] The extension $k(V)/k(X)$ is Galois and its 
Galois group $G$ is cyclic of order $b$. There exists 
$\psi \in k(V)^\times$ such that $\psi^b=\varphi$ and a 
generator of $G$ acts by $\psi \mapsto \zeta \psi$, 
where $\zeta\in k$ is a fixed primitive $b^{th}$-root of 
unity. 

\item[(2)] The relative log pair $h\colon (V,B_V)\to Y$ 
satisfies all properties of a $K$-trivial fibration, 
except that $\rank f_*\calO_X(\lceil \bA(V,B_V)\rceil)$ 
might be bigger than one.

\item[(3)] Both $f\colon (X,B)\to Y$ and 
$h\colon (V,B_V)\to Y$ induce the same discriminant and 
moduli part on $Y$. 

\item[(4)] The group $G$ acts naturally on 
$h_*\calO_V(K_{V/Y})$. The eigensheaf corresponding to the
eigenvalue $\zeta$ is 
$\calL:=f_*\calO_X(\lceil -B+f^*B_Y+f^*M_Y\rceil)\cdot \psi.
$ 
\item[(5)] Assume that $h\colon V\to Y$ is semi-stable 
in codimension one. Then $M_Y$ is an integral divisor,
$\calL$ is semi-positive and 
$
\calL=\calO_Y(M_Y)\cdot \psi.
$
\end{itemize}
\end{lem}

\begin{proof} 
(2) We have $K_V+B_V+(\psi)=h^*(K_Y+B_Y+M_Y)$, and
clearly $(V,B_V)$ has Kawamata log terminal singularities
over the generic point of $Y$.
The generic fibre $H$ of $h$ is a non-singular birational
model of the normalization of $k(F)$ in 
$k(F)((\varphi|_F)^\frac{1}{b})$. Since $b$ is minimal with
$b(K_F+B_F)\sim 0$, $H$ is connected. Therefore 
$\calO_Y=h_*\calO_V$, i.e. $h$ is a contraction.

(3) It follows from (2) and Lemma~\ref{kfbc}. Note that
the assumption $\rank f_*\calO_X(\lceil \bA(V,B_V)\rceil)=1$
is not required in the definition of the discriminant and 
moduli part.

(4) The group $G$ acts on
$\tilde{f}_*\calO_{\tilde{X}}(K_{\tilde{X}/Y})$. Its  
decomposition into eigensheaves is
$$
\tilde{f}_*\calO_{\tilde{X}}(K_{\tilde{X}/Y})=
\bigoplus_{i=0}^{b-1} 
f_*\calO_X(\lceil (1-i)K_{X/Y}-iB+if^*B_Y+if^*M_Y\rceil)
\cdot \psi^i
$$
Since $B-f^*(B_Y+M_Y)$ is supported by the simple normal
crossings divisor $\Sigma_X$, $\tilde{X}$ has rational 
singularities. In particular, 
$h_*\calO_V(K_{V/Y})=\tilde{f}_*\calO_
{\tilde{X}}(K_{\tilde{X}/Y})$ is independent of the choice
of $V$.

(5) By the semi-stable assumption, there exists a big open 
subset $Y^\dagger \subseteq Y$ such that 
$(-B_V+h^*B_Y)|_{h^{-1}(Y^\dagger)}$ is effective
and supports no fibres of $h$.
Since $(\psi|_F)+K_H=-B_H\ge 0$, $\psi$ is a rational section
of $h_*\calO_V(K_{V/Y})$. Furthermore, $\psi \mapsto \zeta \psi$
implies that $\psi$ is a rational section of $\calL$.
Therefore $\calL\subseteq k(Y)\psi$, since $\calL$ has rank 
one by (vi) and (4). We have
$
(h^*a\cdot \psi)+K_{V/Y}=h^*((a)+M_Y)+(-B_V+h^*B_Y).
$

Since $-B_V+h^*B_Y$ is effective over $Y^\dagger$, we infer
that $\calO_Y(M_Y)\psi|_{Y^\dagger} \subseteq 
h_*\calO_V(K_{V/Y})|_{Y^\dagger}$.
Therefore 
$\calO_Y(M_Y)\psi|_{Y^\dagger} \subseteq \calL|_{Y^\dagger}$.
Conversely, let $h^*a\cdot\psi$ be a section of $\calL$. Then
$h^*a\cdot\psi$ is a section of $h_*\calO_V(K_{V/Y})$, i.e.
$(h^*a\cdot \psi)+K_{V/Y} \ge 0$.
Since $-B_V+h^*B_Y$ contains no fibres over codimension one 
points of $Y$, this implies $(a)+M_Y\ge 0$. In particular,
$\calL \subseteq \calO_Y(M_Y)\psi$. Therefore
$
\calO_Y(M_Y)\psi|_{Y^\dagger}= \calL|_{Y^\dagger}.
$
Since $Y^\dagger \subseteq Y$ is a big open subset, this implies
$
\calL^{**}=\calO_Y(M_Y) \psi.
$
By Theorem~\ref{ubc}, $h_*\calO_V(K_{V/Y})$ is locally
free and semi-positive. Its direct summand $\calL$ is 
locally free and semi-positive as well, hence the conclusion.

Finally, over each prime divisor $P$ of $Y$ there exists 
a prime divisor $Q$ of $X$ such that $h(Q)=P$ and
$\mult_Q(-B_V+h^*B_Y)=0$. We infer from (2) that
$\mult_Qh^*(M_Y) \in \Z$. But $\mult_Q h^*(P)=1$, therefore
$M_Y$ is an integral Weil divisor.
\end{proof}

\begin{rem}
Let $\gamma\colon Y'\to Y$ be a generically finite morphism 
from a non-singular quasi-projective variety $Y'$. 
Assume there exists a simple normal crossings divisor 
$\Sigma_{Y'}$ which contains $\gamma^{-1}(\Sigma_Y)$ 
and the locus where $\gamma$ is not \'etale. By base 
change, there exists a commutative diagram
\[ \xymatrix{
& V \ar[dl]_g \ar[dd] &  
& V' \ar[ll]_\nu \ar[dd] \ar[dl]_{g'} \\
X\ar[dr]_f & & X'\ar[dr]_{f'} \ar[ll]_\sigma &  & \\
  & Y              &     & Y' \ar[ll]_\gamma
} \]
such that $(V',B_{V'})-\to (X',B_{X'})\to Y'$
satisfies the same properties (i)-(v).
Here $B_{X'},B_{V'}$ are induced by crepant pull back, 
$\Sigma_{X'}\supseteq \sigma^{-1}(\Sigma_X), \Sigma_{V'}
\supseteq \nu^{-1} (\Sigma_V)$ and 
$\varphi'=\sigma^*\varphi \in k(X')^\times$.
We say that {\em the set-up $(V',B_{V'})-\to (X',B_{X'})\to Y'$ 
is induced by $(V,B_V)\to (X,B) \to Y$ via the base change 
$\gamma\colon Y' \to Y$}.
\end{rem}

\begin{prop}\label{nef}
There exists a finite Galois cover $\tau\colon Y'\to Y$ from
a non-singular variety $Y'$ which admits a simple
normal crossings divisor supporting $\tau^{-1}(\Sigma_Y)$
and the locus where $\tau$ is not \'etale, and such that 
$h'\colon V'\to Y'$ is semi-stable in codimension one for some 
set-up $(V',B_{V'})-\to (X',B_{X'})\to Y'$ induced by 
base change.
\end{prop}

\begin{proof} Let $N$ be the positive integer associated
to $V\to Y$ by Theorem~\ref{ss}. By Theorem~\ref{keycover}, 
there exists a finite Galois cover $\tau\colon Y'\to Y$ such that
$\tau^*(\Sigma_Y)$ is divisible by $N$ and there exists a 
simple normal crossings divisor $\Sigma_{Y'}$ containing
$\tau^{-1}(\Sigma_Y)$ and the locus where $\tau$ is not 
\'etale.
By Theorem~\ref{ss}, there exists an induced set-up
$(V',B_{V'})-\to (X',B_{X'})\to Y'$ induced by base change,
so that $h'\colon V'\to Y'$ is semi-stable in codimension one.
\end{proof}

\begin{prop}\label{kabc} Let $\gamma\colon Y'\to Y$ be a 
generically finite projective morphism from a non-singular
variety $Y'$. Assume there exists a simple normal crossings
divisor $\Sigma_{Y'}$ on $Y'$ which contains 
$\gamma^{-1}(\Sigma_Y)$, and the locus where $\gamma$ is not
\'etale. Let $M_{Y'}$ be the moduli part
of the induced set-up $(V',B_{V'})-\to (X',B_{X'})\to Y'$.
Then $\gamma^*(M_Y)=M_{Y'}$.
\end{prop}

\begin{proof}
{\em Step 1:} Assume that $V/Y$ and $V'/Y'$ are semi-stable 
in codimension one. In particular, $M_Y$ and $M_{Y'}$ are 
integral divisors. Since $h$ is semi-stable in codimension one, 
Theorem~\ref{ubc} implies 
$$
h'_*\calO_{V'}(K_{V'/Y'})  \isoto
\gamma^*(h_*\calO_V(K_{V/Y})).
$$
This isomorphism is natural, hence compatible with the
action of the Galois group $G$. We 
have an induced isomorphism of eigensheaves corresponding to
$\zeta$: $\gamma^*\calO_Y(M_Y)\simeq 
\calO_{Y'}(M_{Y'})$.
Therefore $\gamma^*M_Y-M_{Y'}$ is linearly trivial, and
is exceptional over $Y$. Thus
$\gamma^*M_Y=M_{Y'}$.

{\em Step 2:} By Theorem~\ref{ss} and 
Theorem~\ref{keycover}, we can construct a commutative
diagram 
\[ \xymatrix{
 \bar{Y} \ar[d]_\tau  & \bar{Y}' \ar[d]^{\tau'} 
\ar[l]_{\gamma'} \\
Y &                   Y' \ar[l]_\gamma & 
} \]
as in Remark~\ref{comp}, so that $\bar{V}/\bar{Y}$ is 
semi-stable in codimension one for an induced set-up
$(\bar{V},B_{\bar{V}})-\to (\bar{X},B_{\bar{X}}) \to 
\bar{Y}$.

By Theorem~\ref{ss} and Theorem~\ref{keycover}, we 
replace $\bar{Y}'$ by a finite covering so that 
$\bar{V}'/\bar{Y}'$ is semi-stable in codimension one for
an induced set-up
$(\bar{V}',B_{\bar{V}'})-\to (\bar{X}',B_{\bar{X}'}) \to 
\bar{Y}'$. By Step 1, we have
$M_{\bar{Y}'}={\gamma'}^*(M_{\bar{Y}})$. 
Since $\tau$ and $\tau'$ 
are finite coverings, Lemma~\ref{kfbc} implies 
$\tau^*(M_Y)=M_{\bar{Y}}$ and ${\tau'}^*(M_{Y'})=M_{\bar{Y}'}$.
Therefore ${\tau'}^*(M_{Y'}-\gamma^*(M_Y))=0$, which implies
$M_{Y'}=\gamma^*(M_Y)$.
\end{proof}

\begin{proof}(of Theorem~\ref{main})
Let $f\colon (X,B)\to Y$ be a $K$-trivial fibration 
with $b=b(F,B_F)$ and
$$
K+B+\frac{1}{b}(\varphi)=f^*(K_Y+B_Y+M_Y).
$$
We replace $X$ by a resolution, so that $X$ is 
non-singular, quasi-projective, and $B-f^*(B_Y+M_Y)$ 
is supported by a simple normal crosings divisor $\Sigma_X$. 
Let $V$ be a resolution of the normalization of $X$ in 
$k(X)(\varphi^\frac{1}{b})$ such that $B_V$ has
simple normal crosssings support. We may assume that 
$f,h$ are projective morphisms, after a birational base 
change. Then there exists a closed subvariety 
$\Sigma_f \subsetneq Y$ such that 
$(V,B_V)-\to (X,B)\to Y$ satisfies the assumptions of the 
set-up in the beginning of this section, except that 
$\Sigma_f$ may not be the support of a simple normal 
crosssings divisor. Let $\sigma\colon Y'\to Y$ be an 
embedded resolution so that 
$\Sigma_{Y'}:=\sigma^{-1}(\Sigma_f)$ is a divisor with 
simple normal crossings. There exists an induced set-up
$(V',B_{V'})-\to (X',B')\to Y'$. 

We claim that $\sigma^*(\bM_{Y'})=\bM_{Y''}$ and
$\sigma^*(K_{Y'}+\bB_{Y'})=K_{Y''}+\bB_{Y''}$ for 
every birational contraction $\sigma\colon  Y'' \to Y'$.
By Hironaka's resolution of singularities, there
exists a diagram of birational morphisms
\[ \xymatrix{
Y'' \ar[d]_\sigma & Y''' \ar[l] \ar[dl]^{\sigma'} \\
Y' &
} \]
such that $Y'''$ is a non-singular quasi-projective variety
admitting a simple normal crossings divisor which
supports ${\sigma'}^{-1}(\Sigma_{Y'})$ and the exceptional 
locus of $Y'''/Y'$. By Proposition~\ref{kabc}, 
${\sigma'}^*(\bM_{Y'})=\bM_{Y'''}$ and consequently
${\sigma'}^*(K_{Y'}+\bB_{Y'})=K_{Y'''}+\bB_{Y'''}$.
Since $Y'''/Y''$ is a birational morphism, the claim
follows.

Let $\tau\colon \bar{Y}'\to Y'$ be a covering given by
Proposition~\ref{nef}. By Lemma~\ref{ex}, 
$M_{\bar{Y}'}$ is a Cartier divisor and $\calO_{\bar{Y}'}
(M_{\bar{Y}'})$ is a semi-positive invertible sheaf. 
In particular, $M_{\bar{Y}'}$ is nef/$S$. 
But $\tau^*(M_{Y'})=M_{\bar{Y}'}$
according to Lemma~\ref{kfbc}, hence $M_{Y'}$ is nef/$S$.
\end{proof}

\begin{proof}(of Theorem~\ref{cv}) By the Lefschetz principle,
we may assume $k=\mathbb C$. After a finite base change
(Lemma~\ref{kfbc}), we may assume that the induced root fiber 
space $h\colon V\to Y$ is {\em semi-stable}. By construction, 
the invertible sheaf
$\calL:=\calO_Y(M_Y)\subset h_*\omega_{V/Y}$ is a 
direct summand. 

We know that $M_Y$ is a nef Cartier divisor
on the curve $Y$. If $\deg(M_Y)>0$, then $M_Y$ is ample, 
in particular semi-ample. If $\deg(M_Y)=0$, Theorem~\ref{fk}
implies $\calL^{\otimes m}\simeq \calO_Y$, so 
$M_Y$ is semi-ample.
\end{proof}


\section{Asymptotically saturated algebras}


We first recall some terminology from \cite{Plflips}.
Let $\pi\colon X\to S$ be a proper morphism. A {\em normal 
functional algebra} of $X/S$ is an $\calO_S$-algebra 
of the form
$$
\calL=\calR_{X/S}(\bM_\bullet)=
\bigoplus_{i=0}^\infty \pi_*\calO_X(\bM_i),
$$
where $\{\bM_i\}$ is a sequence of b-free/$S$ b-divisors
of $X$ such that $\bM_i+\bM_j\le \bM_{i+j}$ for every $i$
and $j$. The sequence of $\Q$-b-divisors $\bD_i=\frac{1}{i}
\bM_i$ is called the {\em characteristic sequence} of
$\calL$. The algebra $\calL$ is {\em bounded} if there
exists an $\R$-b-divisor $\bD$ of $X$ such that $\bD_i\le \bD$
for every $i$. By the Limiting Criteria~\cite[Theorem 4.28]
{Plflips}, the $\calO_S$-algebra $\calL$ is finitely 
generated if and only if the characteristic sequence 
$\bD_\bullet$ is constant up to a truncation.
For an $\R$-b-divisor $\bA$, the algebra $\calL$ is 
{\em asymptotically $\bA$-saturated}, if there exists a 
positive integer $I$ such that
$$
\pi_*\calO_X(\lceil \bA+j\bD_i\rceil) \subseteq
\pi_*\calO_X(\bM_j)
\mbox{ for } I|i,j.
$$

The {\em Kodaira dimension} of $\calL$ is 
$\kappa(\calL):=\max_i \kappa(X/S,\bM_i)$. 
 We say that $\calL$ is a {\em big algebra} if 
$\kappa(\calL)=\dim(X/S)$.

\begin{defn} Let $\calL$ be a normal functional algebra
of $X/S$. There exists a unique rational map with 
connected fibers $f\colon X -\to Y/S$ and a normal functional 
algebra $\calL'$ of $Y/S$, such that $f^*\colon \calL'\to \calL$ 
is a quasi-isomorphism and $\calL'$ is a big algebra.
We say that $(f,\calL')$ is the {\em Iitaka 
fibration} of $\calL$.
\end{defn}

\begin{proof}\cite[Lemma 6.22]{Plflips}
Let $\calL=\calR_{X/S}(\bM_\bullet)$.
Since $\bM_i+\bM_j\le \bM_{i+j}$ and the $\bM_i$'s
are b-free, there exists $I \in \N$ and
a rational map $f\colon X -\to Y/S$ which is the Iitaka 
contraction of $\bM_i$ for every $i$ divisible by $I$. 
Up to a quasi-isomorphism, we may assume that the b-free 
b-divisors $\bM_i$ are effective. Since $f$ has connected
fibers, there exists a convex sequence $\bM'_\bullet$ such 
that $\bM_i=f^*(\bM'_i)$ for every $I|i$.
In particular, $\calL$ is quasi-isomorphic to the big algebra
$\calL':=\calR_{Y/S}(\bM'_\bullet)$.
\end{proof}

\begin{lem}\label{rkone} 
Let $(f\colon X -\to Y/S,\calL')$ be the Iitaka
fibration of a normal functional algebra $\calL$.
If $\calL$ is asymptotically $\bA$-saturated, then 
$\rank f'_*\calO_{X'}(\lceil \bA \rceil )\le 1$, 
where $f'\colon X'\to Y$ is a regular representative of 
the rational function $f$.
\end{lem}

\begin{proof} We may assume that $f'=f$ and 
$\calL=\calR_{X/S}(f^*\bM'_\bullet)$, where
$\calL'=\calR_{Y/S}(\bM'_\bullet)$ is the induced
big algebra. By assumption, there exists $i$ such 
that $\bD_i$ is b-big/$S$. After passing to higher 
models, we may assume that $\bD_i$ descends to $Y$. 
There exists a birational contraction $\mu\colon Y\to Z/S$ 
and an ample/$S$ $\Q$-divisor $H$ on $Z$ such that
$(\bD_i)_Y\sim_\Q \mu^*H$. For $j$ sufficiently large
and divisible, the $\calO_Z$-sheaf
$$
\mu_*(f_*\calO_X(\lceil \bA +jf^*\bD_i\rceil))=
\mu_*f_*\calO_X(\lceil \bA \rceil)\otimes \calO_Z(jH)
$$
is $\pi$-generated. Therefore
$f_*\calO_X(\lceil \bA +jf^*\bD_i\rceil)$ is 
generically $\pi$-generated. Asymptotic saturation 
implies that 
$f_*\calO_X(\lceil \bA +jf^*\bD_i\rceil)$ is 
contained in the b-divisorial sheaf $\calO_Y(\bM_j)$
on an open subset of $Y$. The latter has rank 
one, hence $f_*\calO_X(\lceil \bA \rceil)$ has rank
at most one.
\end{proof}

\begin{prop}\label{div} 
(cf. ~\cite[Proposition 4.50]{Plflips})
Consider a commutative diagram
\[ \xymatrix{
(X,B) \ar[dr]_\pi \ar[rr]^f & & Y \ar[dl]^\sigma \\
& S  &
} \]
and a normal functional algebra 
$\calL=\calR_{X/S}(\bM_\bullet)$ with the following 
properties:
\begin{itemize}
\item[(a)] $f\colon (X,B)\to Y$ is a $K$-trivial 
fibration.
\item[(b)] $\calL$ is bounded and asymptotically
$\bA(X,B)$-saturated.
\item[(c)] There exist b-divisors $\bM'_i$ of $Y$ such 
that $\bM_i=f^*(\bM'_i)$ for all $i$.
\end{itemize}

Then $\calL':=\calR_{Y/S}(\bM'_\bullet)$ is a normal 
bounded functional algebra of $Y/S$, which is asymptotically
$\bA_{div}$-saturated. Moreover, the natural map 
$f^*\colon \calL'\to \calL$ is an isomorphism of 
$\calO_S$-algebras.
\end{prop}

\begin{proof} It is clear that $\calL'$ is a functional
algebra of $Y/S$, and $f^*\colon \calL'\to \calL$ is an
isomorphism of $\calO_S$-algebras. The algebra is normal
since each $\bM'_i$ is b-free. Let $\bD'_i=\frac{1}{i}\bM'_i$ 
be the characteristic sequence of $\calL'$. 

We first check that $\calL'$ is bounded.
After passing to higher models, we may assume that 
there exists an effective Cartier divisor 
$E$ on $X$ such that $f^*\bD'_i=\bD_i \le \overline{E}$.
Let $E'$ be the divisorial support of $f(\Supp(E^v)) \subset Y$.
For each $i$, we can find a birational model $X'/Y'$ of 
$X/Y$, fitting in the commutative diagram
\[ \xymatrix{
X \ar[d]_f  & X' \ar[d]^{f'} \ar[l]_h \\
Y           &   Y'  \ar[l]  
   } \]
such that $\bD'_i$ and $\bD_i$ descend on $Y'$ and $X'$
respectively. In particular, ${f'}^*((\bD'_i)_{Y'}) \le h^*E$. 
Since $\bD'_i$ is effective and $Y'/Y$ is an isomorphism over 
a big open subset of $Y$, we conclude that $(\bD'_i)_Y$ is 
supported by $E'$. This holds for every $i$, hence $\bD'_\bullet$ 
is bounded.

It remains to check asymptotic $\bA_{div}$-saturation.
Fix two integers $i,j$ which are divisible by $I$. 
By Theorem~\ref{main}, we may assume the following
properties hold (after a birational base change):
\begin{itemize}
\item[(i)] $X,Y$ are non-singular.
\item[(ii)] $\bD'_i,\bD'_j$ descend to $Y$ (in particular
$\bD_i,\bD_j$ descend to $X$). Denote $D'_i=(\bD'_i)_Y$
and $D'_j=(\bD'_j)_Y$.
\item[(iii)] $\Supp(B)\cup \Supp(f^*D'_i)$ and
$\Supp(\bB_Y)\cup \Supp(D'_i)$ are simple normal crossings
divisors on $X$ and $Y$, respectively.
\item[(iv)] $\bA_{div}=\bA(Y,\bB_Y)$.
\end{itemize}
Under these assumptions, the saturation for $i,j$ means
$$
\pi_*\calO_X(\lceil -B + jf^*D'_i \rceil)
\subseteq \pi_*\calO_X(jf^*D'_j)
$$ 
By Lemma~\ref{a1} and Remark~\ref{a2}, 
$$
f_*\calO_X(\lceil -B + jf^*D'_i \rceil)=
\calO_Y(\lceil -\bB_Y+jD'_i\rceil).
$$
Since $\pi_*\calO_X(jf^*D'_j)=\sigma_*\calO_Y(\bM'_j)$,
we infer
$$
\sigma_*\calO_Y(\lceil -\bB_Y + jD'_i \rceil)
\subseteq \sigma_*\calO_Y(jD'_j).
$$
Therefore
$
\sigma_*\calO_Y(\lceil \bA_{div}+j\bD'_i\rceil)
\subseteq 
\sigma_*\calO_{Y'}(\bM'_j),
$
by (i)-(iv) again.
\end{proof}

\begin{exmp} (Reduction to big algebras) Let
$(X/S,B)$ be a relative log pair, and let $\calL$ be 
a normal, bounded functional algebra with
Iitaka fibration $(f\colon X-\to Y/S,\calL')$, satisfying 
the following properties:
\begin{itemize}
\item[(i)] $K_{X'}+B_{X'}\sim_\Q {f'}^*D$, where 
$f'\colon X'\to Y$ is a regular model of $f$ and 
$B_{X'}$ is a crepant boundary ($\bA(X,B)=\bA(X',B_{X'})$).
\item[(ii)] $(X',B_{X'})$ has klt singularities over 
the generic point of $Y$. 
\item[(iii)] $\calL$ is asymptotically 
$\bA(X,B)$-saturated.
\end{itemize}

 By Lemma~\ref{rkone}, $f'\colon (X',B_{X'}) \to Y$
is a $K$-trivial fibration. By Theorem~\ref{main} and 
Proposition~\ref{div}, we may replace $Y$ by a higher 
birational model so that the following properties hold:

\begin{itemize}
\item[(a)] $\bA_{div}=\bA(Y,\bB_Y)$.

\item[(b)] $K_{X'}+B_{X'}\sim_\Q {f'}^*(K_Y+\bB_Y+\bM_Y)$.

\item[(c)] $\calL'$ is normal, bounded 
and asymptotically $\bA(Y,\bB_Y)$-saturated. 
\end{itemize}
\end{exmp}

The above example is a first step towards a reduction of 
(0LP) algebras~\cite[Remark 4.40]{Plflips}) to the big 
case. To complete the reduction, we need to know that the 
moduli b-divisor $\bM$ is b-semi-ample.
However, the b-nef property of the moduli b-divisor
is enough for some applications to the Fano case.
We show that the restriction of an FGA algebra to an
exceptional log canonical centre is again an FGA 
algebra (cf. ~\cite[Proposition 4.50]{Plflips} for 
lc centers of codimension one):

\begin{thm}\label{dj} 
Let $(X/S,B)$ be a relative generalized Fano log variety, 
let $\nu\colon W\to X$ be the normalization of an exceptional 
lc centre of $(X,B)$ and let $\calL=\calR_{X/S}(\bM_\bullet)$ 
be a normal bounded functional algebra of $X/S$ such that 
the following hold:

\begin{itemize}
\item[(i)] $\calL$ is asymptotically 
$(\bA(X,B)+E)$-saturated, where $E$ is the unique lc place
over $W$.
\item[(ii)] There exists an open subset $U\subseteq X$ such 
that $U\cap \nu(W)\ne \emptyset$, 
$\bD_i|_U=\overline{D}|_U \ \ \forall i$ for some $\Q$-Cartier 
divisor $D$ on $X$.

\item[(ii$^*$)] $U$ contains $(X,B)_{-\infty}\cap \nu(W)$ 
and $C \cap \nu(W)$, for every lc centre $C\ne \nu(W)$ of 
$(X,B)$.
\end{itemize}

Then there exists a well defined {\em restricted algebra}
$\calL\frest W$ of $W/S$, with the following properties:

\begin{itemize}
\item[(1)] $\calL\frest W=\calR_{W/S}(\bM'_i)$ is a normal, 
bounded functional algebra.

\item[(2)] $\calL\frest W$ is $\bA(W,B_W)$-saturated, where
$(W/S,B_W)$ is a relative generalized Fano log variety.

\item[(3)] $\LCS(W,B_W)\subset U':=U|_W$ and 
$\bD'_i|_{U'}=\overline{D|_{U'}}$ for every $i$. 

\item[(4)] The $\calO_S$-algebras $\calL\frest W$ and 
$\calL\frest E$ are quasi-isomorphic. 

\end{itemize}
\end{thm}

\begin{proof} Let $H$ be an ample/$S$ $\Q$-divisor on $X$
such that $-(K+B+H)$ is ample/$S$. Then 
~\cite[Theorem 4.9]{qlv} constructs an effective $\Q$-divisor 
$B_W$ on $W$ such that $(W/S,B_W)$ is a relative generalized 
Fano log variety with $(K+B+H)|_W\sim_\Q K_W+B_W$ and 
$\LCS(W,B_W)$ is contained in the union of $(X,B)_{-\infty}$ 
and all lc centres of $(X,B)$ different than $\nu(W)$. In 
particular, $\LCS(W,B_W)\subset U'$. Consider the induced 
diagram:
$$
{{
\begin{array}{ccc}
(E,B_E) &\subset &(X',B_{X'}) \\
\downarrow & &\downarrow \\
W &\to & (X,B)
\end{array}
}}
$$
By adjunction and Kawamata-Viehweg vanishing, $\calL\frest E$
is asymptotically $\bA(E,B_E)$-saturated~\cite[Proposition 4.50]
{Plflips}.
By (ii), there exist b-free/$S$ b-divisors $\bM'_i$ of $W$
such that $\bM_i\frest E=h^*(\bM'_i)$ for every $i$.
By construction, $(E,B_E)\to W/S$ is a $K$-trivial fibration
for which Proposition~\ref{div} applies. Therefore 
$\calL\frest W:=\calR_{W/S}(\bM'_i)$ is quasi-isomorphic to 
$\calL\frest E$, it is normal, bounded and asymptotically 
$\bA_{div}$-saturated. 
From the construction of $B_W$ (choosing $W'$ high enough so 
that $\bA_{div}=\bA(W',\bB_{W'})$, in the proof of 
~\cite[Theorem 4.9]{qlv}), we have $\bA(W,B_W)\le \bA_{div}$.
Therefore $\calL\frest W$ is asymptotically $\bA(W,B_W)$-saturated.

Finally, $\bD_i|_U=\overline{D|_U}$ implies
$\bD'_i|_{U'}=\overline{D|_{U'}}$.
\end{proof}


\section{Parabolic fiber spaces}


A {\em parabolic fiber space} is a contraction
of non-singular proper varieties $f\colon X\to Y$ such that
the generic fiber $F$ has Kodaira dimension zero.
Let $b$ be the smallest positive integer with
$|bK_F|\ne \emptyset$. We fix a rational function 
$\varphi \in k(X)^\times$ such that 
$K+\frac{1}{b}(\varphi)$ is effective over the generic 
point of $Y$.

\begin{defnprop}(cf.~\cite[Corollary 2.5]{fm}) \label{pb}
Let $f\colon X\to Y$ be a parabolic fiber space with a choice 
of a rational function $\varphi$, as above. There exists 
a unique b-nef $\Q$-b-divisor $\bM=\bM(f,\varphi)$ of 
$Y$ satisfying the following properties:
\begin{itemize}
\item[(1)] Let $\varrho\colon Y'\to Y$ be a surjective proper
morphism, and let $f'\colon X'\to Y'$ be an induced parabolic 
fiber space:
\[ 
\xymatrix{
 X \ar[d]_f  & X' \ar[l]_\nu    \ar[d]_{f'}   \\
 Y  & Y' \ar[l]_\varrho
} \]
Then $\varrho^*\bM(f,\varphi)\sim_\Q \bM(f',\nu^*\varphi)$.
Moreover, $\varrho^*\bM(f,\varphi)=\bM(f',\nu^*\varphi)$
if $\varrho$ is generically finite.

\item[(2)] If $f$ is semi-stable in codimension one, then
$$
\bigoplus_{b|i} f_*\calO_X(iK_X)^{**} =
\bigoplus_{b|i} \calO_Y(i(K_Y+\bM_Y)) \cdot \varphi^i.
$$
\end{itemize}
We say that $\bM=\bM(f,\varphi)$ is the 
{\em moduli $\Q$-b-divisor} associated to the parabolic 
fiber space $f$. If $\varphi'$ is another choice
of the rational function, then $b\bM\sim b\bM'$.
Therefore $b\bM$ is uniquely defined up to linear
equivalence. 
\end{defnprop}

\begin{proof} 
There exists~\cite[Proposition 2.2]{fm} a unique 
$\Q$-divisor $B_X$ on $X$ satisfying the following 
properties:
\begin{itemize}
\item[(i)] $K+B_X+\frac{1}{b}(\varphi)=f^*D$ for some 
$D\in \Div(Y)_\Q$.
\item[(ii)] There exists a big open subset 
$Y^\dagger \subseteq Y$ such that $-B_X|_{f^{-1}(Y^\dagger)}$ 
is effective and contains no fibers of $f$ in its support.
\end{itemize}

One can easily check that $f\colon (X,B_X)\to Y$ is 
a $K$-trivial fibration.
Let $\bB, \bM$ be the induced discriminant and moduli 
b-divisors of $Y$. Note that if $X-\to X'$ is a birational 
map over $Y$, then $(X,B_X)-\to (X',B_{X'})$ is crepant, 
i.e. $\bA(X,B)=\bA(X',B_{X'})$. Also, $\lfloor \bB 
\rfloor=0$. We claim that $\bM=\bM(f,\varphi)$ (the 
uniqueness is clear by semi-stable reduction in 
codimension one). Since $\kappa(F)=0$, and
$$
K+B_X+ \frac{1}{b}(\varphi)=f^*(K_Y+\bB_Y+\bM_Y),
$$
we infer that
$
f_*\calO_X(iK)^{**}=
\calO_Y(i(K_Y+\bB_Y+\bM_Y))\cdot \varphi^i \mbox{ for } b|i.
$
If $f$ is semi-stable in codimension one, then 
$\bB_Y=0$. Therefore (2) holds.

By Theorem~\ref{main}, $\bM$ is b-nef. 
We may replace $f$ by a birational base change so
that $\bM$ descends to $Y$: $\bM=\overline{\bM_Y}$.
If $\varrho$ is generically finite, (1) holds by 
inversion of adjunction for finite morphisms 
(cf.~\cite[Proposition 4.7]{fm}). For general $\varrho$,
we imitate the proof of Proposition~\ref{kabc} with the
following simplifications: the root fiber space $h\colon V\to Y$ 
is parabolic as well; if $h$ has simple normal 
crossings degeneration and is semi-stable in codimension
one, then $\bM_Y$ is an integral divisor and 
$
h_*\calO_V(K_{V/Y})=\calO_Y(\bM_Y)\cdot \varphi.
$
\end{proof}

A non-singular projective variety $X$ 
{\em has a good minimal model} if it is birational
to a normal projective variety $Y$ such that
$Y$ has terminal singularities and $K_Y$ is semi-ample.
This holds if $\dim(X)\le 3$, by the Minimal Model 
Program and Abundance (see ~\cite{mmp}). 
The following is a generalization 
of ~\cite[Theorem 13]{koddim}:

\begin{thm} \label{k0}
Let $f\colon X\to Y$ be a parabolic fiber space such that
its geometric generic fibre 
$\bar{F}=X\times_Y \Spec(\overline{k(Y)})$ has
a good minimal model over $\overline{k(Y)}$. 
Then there exists a diagram
\[ \xymatrix{
X \ar[d]_f & \bar{X}\ar[l]\ar[d]_{\bar{f}}
\ar[r] & X^!\ar[d]_{f^!}\\
Y        & \bar{Y} \ar[l]_\tau \ar[r]^\varrho & Y^!
} \]
such that the folowing hold:
\begin{itemize}
\item[(1)] $\bar{f}$ and $f^!$ are parabolic fiber spaces.
\item[(2)] $\tau$ is generically finite, and $\varrho$
is a proper dominant morphism.
\item[(3)] $\bar{f}$ is birationally induced via base 
change by both $f$ and $f^!$.
\item[(4)] $\tau^*\bM=\bar{\bM}\sim_\Q \varrho^*\bM^!$,
where $\bM,\bar{\bM},\bM^!$ are the corresponding moduli 
$\Q$-b-divisors.
\item[(5)] $\bM^!$ is b-nef and big and 
$\Var(f^!)=\dim(Y^!)$.
\end{itemize}

In particular, $\kappa(\bM)=\Var(f)$, where
$\Var(f)$ is the variation of the fiber space $f$.
\end{thm}

\begin{proof} By the definiton of the variation of a 
fibre space, there exists a diagram as above, satisfying 
(1), (2),(3) and such that $\dim(Y^!)=\Var(f)=\Var(f^!)$. 
After a generically finite base change, we may also assume 
that $\bM^!$ descends to $Y^!$, and $f^!$ is semi-stable in 
codimension one. 

By (3) and Definition-Proposition~\ref{pb}, (4) holds. In 
particular, $\kappa(\bM)=\kappa(\bM^!)$. Since 
$\bar{F}$ has a good minimal model, Viehweg's $Q(f^!)$ 
Conjecture holds~\cite[Theorem 1.1.(i)]{minit}, that is 
the sheaf
$
(f^!_*\omega_{X^!/Y^!}^i)^{**}
$
is big for $i$ large and divisible. But 
$
(f^!_*\omega_{X^!/Y^!}^i)^{**} \simeq \calO_{Y^!}
(i\bM^!_{Y^!})
$
for $b|i$, since $f^!$ is semi-stable in codimension 
one. Equivalently, $\kappa(Y^!,\bM^!_{Y^!})=
\dim(Y^!)$. Therefore $\bM^!$ is b-nef and big.
\end{proof}


\end{document}